\documentclass[12pt]{elsart}

\usepackage{amsfonts}
\usepackage{amsmath}
\usepackage{amssymb}

\numberwithin{equation}{section}

\renewcommand\thefigure{\thesection.\@arabic\c@figure}
\renewcommand\thetable{\thesection.\@arabic\c@table}

\def\eps{\varepsilon}

\def\s{\sigma}
\def\qed{\hfill\rule{.2cm}{.2cm}}
\def\square{\ifmmode\sqr\else{$\sqr$}\fi}
\def\sqr{\vcenter{
         \hrule height.1mm
         \hbox{\vrule width.1mm height2.2mm\kern2.18mm\vrule width.1mm}
         \hrule height.1mm}}                  
\def\proof{\noindent{\bf Proof. }}

\def\P{{\mathbb P}}
\def\E{{\mathbb E}}
\def\Z{{\mathbb Z}}
\def\N{{\mathbb N}}
\def\R{{\mathbb R}}

\def\p{{\mathcal P}}

\def\1{{\bf 1}}

\def\a{{\alpha}}
\def\b{{\beta}}

\def\ln{{\l_n}}
\def\bln{{\bar\l_n}}
\def\l{{\lambda}}
\def\tl{{\tilde\lambda}}

\def\reff#1{(\ref{#1})}

\newtheorem{theo}{Theorem}[section]
\newtheorem{corol}[theo]{Corollary}
\newtheorem{lm}[theo]{Lemma}

\newtheorem{rmk}[theo]{Remark}

\begin{document}
\bibliographystyle{alpha}
\runauthor{Ferrari, Fontes, Niederhauser and Vachkovskaia}

\begin{frontmatter}
\title{The serial harness interacting with a wall}
\date{June 2, 2003}

\author[USP]{Pablo A.~Ferrari}
\author[USP]{Luiz R.~G.~Fontes}
\author[USP]{Beat M.~Niederhauser}
\author[UNICAMP]{Marina Vachkovskaia}
\address[USP]{Universidade de S\~{a}o Paulo}
\address[UNICAMP]{Universidade Estadual de Campinas}



\begin{abstract}
The serial harnesses introduced by Hammersley describe the motion of a
hypersurface of dimension $d$ embedded in a space of dimension $d+1$. The
height assigned to each site $i$ of $\Z^d$ is updated by taking a weighted
average of the heights of some of the neighbors of $i$ plus a ``noise'' (a
centered random variable). The surface interacts by exclusion with a ``wall''
located at level zero: the updated heights are not allowed to go below zero.
We show that for any distribution of the noise variables and in all
dimensions, the surface delocalizes. This phenomenon is related to the so
called ``entropic repulsion''. For some classes of noise distributions,
characterized by their tail, we give explicit bounds on the speed of the
repulsion.
\end{abstract}



\end{frontmatter}


\section{Introduction and results}
Hammersley (1965) introduced the \emph{serial harness}, a discrete-time
stochastic process that models the time evolution of a hypersurface of
dimension $d$ embedded in a $d+1$ dimensional space. A quantity $Y_n(i)\in \R$
stays for the height of the surface at site $i\in \Z^d$ at (integer) time
$n\ge 0$. The initial configuration is the flat surface $Y_0(i)=0$ for all
$i$.  Under the evolution, at each moment $n\ge 0$ the height at each site is
substituted by a weighted average of the heights at the previous moment plus a
symmetric random variable.

Let ${\mathcal P}=\{p(i,j)\}_{i,j \in \Z^d}$ be a stochastic matrix, i.e.\
$p(i,j) \geq 0$ and $\sum_j p(i,j) = 1$, which satisfies $p(i,j) = p(0,
j-i)=:p(j-i)$ (homogeneity), $\sum_j j p(j) = 0$, and $p(j) = 0$ for all $|j|
> v$ for some $v$ (finite range). Assume also that ${\mathcal P}$ is truly
$d$-dimensional: $\{j\in \Z^d\,:\, p(j)\neq 0\}$ generates $\Z^d$.

Let ${\mathcal E}=(\eps, (\eps_n(i),\, i \in \Z^d),n \in \Z)$ be a
family of i.i.d.\  integrable symmetric
random variables.
Let $\P$ and $\E$ denote the probability and expectation in the
probability space generated by ${\mathcal E}$. (We use preliminary $n\in \N$ in the
definitions but later it will be useful to have $n\in \Z$.)

The \emph{serial harness} $(Y_n,\, n\ge 0)$ is the discrete-time Markov
process in $\R^{\Z^d}$ defined by
\begin{eqnarray}
\label{I.21}
Y_n(i) = \left\{
           \begin{array}{lll}{}\hspace{2cm} 0, & \mbox{ if }& n =0, \\
               \sum_{j \in \Z^d} p(i,j) Y_{n-1}(j) +  \eps_n(i),
                        & \mbox{ if } & n \geq 1. \\
         \end{array}
         \right.
\end{eqnarray}
Here $Y_n(i)$ denotes the height of the serial harness at site $i$ at time
$n$.  In other words, the evolution is given by
\begin{equation}
  \label{I.23a}
Y_n={\mathcal P}Y_{n-1}+\eps_n.
\end{equation}
where $\eps_n = (\eps_n(i)\,,\, i\in \Z^d)$. Since the ``noise variable''
$\eps$ is symmetric and thus has zero mean, we have that $\E Y_n(i) = 0$ for
all $i,n$.  We can interpret $p(i,j)$ as transition probabilities of a random
walk on $\Z^d$; let $p_m(i,j)$ be its $m$-step transition probabilities. By
homogeneity, $p_m(i,j)=p_m(0,j-i)=:p_m(j-i)$. Iterating \reff{I.21},
\begin{equation}
\label{I.23}
Y_n(i)\; =\; \sum_{r=1}^n \sum_{j \in \Z^d}p_{n-r}(i,j) \eps_r(j)
        \;\stackrel{d}=\;\sum_{r=0}^{n-1}\sum_{j\in\Z^d}p_{r}(j)\eps_r(j),
\end{equation}
for all $n\geq1,i\in\Z^d$, where $\stackrel{d}=$ means equidistributed.
Hammersley (1965) obtained that
\begin{equation}
  \label{ham1}
 \E (Y_n(i))^2 = \sigma^2 s(n)
\end{equation}
where $\sigma^2$ is the variance of $\eps$ and
\begin{equation}
\label{I.23bis}
s(n) := \sum_{r=0}^{n-1} \sum_{j \in \Z^d} p_r(j)^2.
\end{equation}
is the expected number of encounters up to time $n$ of two independent copies
of a random walk starting at 0 with transition probabilities ${\mathcal P}$.
Equality \reff{ham1} follows immediately from \reff{I.23}. Since $s(n)\sim
\sqrt n$ for $d=1$, $s(n) \sim \log n$ for $d=2$ and $s(n)$ is uniformly
bounded in $n$ for $d\ge 3$ (see, for example, Spitzer (1976)), the surface
delocalizes in dimensions $d\le 2$ and stays localized in dimensions $d\ge 3$.
Toom (1997) studies localization of the surface and surface-differences in
function of the decay of the distribution of $\eps$.

We consider the serial harness interacting by exclusion with a wall located at
the origin. The \emph{wall process} $(W_n,\,n\ge 0)$ is the Markov process in
$(\R^+)^{\Z^d}$ defined by
\begin{eqnarray}
\label{I.22}
W_n(i) = \left\{
            \begin{array}{lll}{}\hspace{2cm}  0, & \mbox{ if } & n =0, \\
                \Big(\sum_{j \in \Z^d}
                        p(i,j) W_{n-1}(j)\, +\, \eps_n(i)\Big)^+,
                        & \mbox{ if }  & n \geq 1, \\
         \end{array}
         \right.
\end{eqnarray}
for $i\in \Z^d$, where for $a\in\R$, $a^+=a\vee0=\max(a,0)$; this can be
reexpressed as
\begin{equation}
  \label{I.22a}
W_n=({\mathcal P}W_{n-1}+\eps_n)^+.
\end{equation}
We say that the law of a random surface $Z$ is an \emph{invariant measure} for
the wall process if $Z\stackrel{d}=(\eps_0+\p Z)^+$, with $\eps_0$ and $Z$
independent. We show in Section \ref{s2} that
\begin{equation}
  \label{p11}
  W_n\le W_{n+1} \hbox{ stochastically}.
\end{equation}
This implies that $W_n$ is stochastically non-decreasing and thus their laws
converge to a limit (that could give positive weight to infinity). If the
limit is nondegenerate, then it is an invariant measure for the wall process.
Monotonicity \reff{p11} implies in particular \[\mu_n:=\E W_n(0)\] is
nondecreasing and thus converges either to a finite limit or to $\infty$.  Our
first result is general and rules out the former possibility, showing however
that $\mu_n$ goes to infinity slower than $n$.

\begin{theo}
\label{delocalisation}
(a) There is no nondegenerate invariant measure for the wall process
  $(W_n)$;
(b) $W_n\to\infty$ in probability;  (c) $\mu_n\to\infty$ as $n\to\infty$; (d)
  $\mu_n/n \to 0$ as $n\to\infty$.
\end{theo}

This theorem is proven in Section \ref{s2}.

Let $F$ be the law of
$\eps$, $\bar F(x) = \P(\eps >x)$ and define
\begin{eqnarray}
  {\mathcal L}_\a^-&:=&\{F:\bar F(x)\leq ce^{-c'x^{\a}}, x>0,
                     \mbox{ for some positive $c$, $c'$}\}\label{eq:logpol-}\\
 {\mathcal L}_\a^+&:=&\{F:\bar F(x)\geq ce^{-c'x^{\a}}, x>0,
                     \mbox{ for some positive $c$, $c'$}\}\label{eq:logpol+}
\end{eqnarray}
and
\begin{eqnarray}
{\mathcal L}_\a&:=&{\mathcal L}_\a^-\cap{\mathcal L}_\a^+ \label{eq:logpol}
\end{eqnarray}

We next state our main result. It consists of upper and lower bounds for
$\mu_n$ for different noise distributions.
\begin{theo}
\label{bounds_gauss}
There exist constants $c$ and $C$ that
may depend on the dimension
such that
\begin{itemize}
\item[(i)] for $d=1$ if $F\in {\mathcal L}_1^-$,
  \begin{equation}
    \label{eq:b1}
c n^{1/4} \leq \mu_n \leq C n^{1/4}\sqrt{\log n};
  \end{equation}
\item[(ii)] for $d=2$, if $F\in {\mathcal L}_\a$, for some $\a\ge 1$,
\begin{equation}
    \label{eq:b2}
    c (\log n)^{\frac1\a\vee\frac12}\leq \mu_n \leq C \log n;
  \end{equation}
\item[(iii)] for  $d\geq3$, if $F\in {\mathcal L}_\a$, for some $1\le\a\ne1+d/2$,
\begin{equation}
    \label{eq:b3}
    c (\log n)^{\frac1\a}\leq\mu_n\leq C(\log n)^{\frac1\a\vee\frac2{2+d}};
  \end{equation}
\item[(iv)] for  $d\geq3$ if $F\in {\mathcal L}_{1+d/2}$,
\begin{equation}
    \label{eq:b4}
   c (\log n)^{\frac{2}{2+d}}\leq\mu_n \leq C (\log
   n)^{\frac{2}{2+d}}(\log\log n)^{\frac{d}{2+d}}.
  \end{equation}
\end{itemize}
\end{theo}

Our upper bound in \reff{eq:b4} can be slightly improved, see
\reff{eq:ub2} and Remark~\ref{rmk:ub2} below.
The lower bound in (i) can be shown to hold under
weaker conditions; that is also the case for some cases of (ii);
see~\reff{I.24} and Remark~\ref{4moma} below.
If the noise distribution is in ${\mathcal L}_\a$ for some
$\a\geq1$, then our lower and upper bounds to $\mu_n$ are of the same
order in the case that $d\ge 3,\,1\le \alpha< 1+d/2$ (which includes the
Gaussian case $\alpha=2$ for all such dimensions), and also in the case that
$d=2,\,\alpha=1$.

Theorems \ref{delocalisation} and \ref{bounds_gauss} catch the effect of the
``entropic repulsion'' in a stochastically moving surface interacting with a
wall by exclusion.

Many papers deal with the problem of entropic repulsion in Equilibrium
Statistical Mechanics. The role of the entropic repulsion in the Gaussian free
field was studied by Lebowitz and Maes (1987), Bolthausen, Deuschel and
Zeitouni (1995), Deuschel (1996), Deuschel and Giacomin (1999) and
Bolthausen, Deuschel and Giacomin (2001).  In the Ising, SOS and related
models the matter was discussed in Bricmont, El Mellouki and Fr\"ohlich
(1986), Bricmont (1990), Cesi and Martinelli (1996), Dinaburg and Mazel (1994), Holick\'y and Zahradn\'\i k (1993),
and Ferrari and Mart\'{\i}nez (1998).

The exponent $1/4$ for dynamic entropic repulsion in $d=1$ was predicted by
Lipowsky (1985) using scaling arguments. This exponent was then found
numerically by Mon, Binder, Landau (1987), Binder (1990), De Coninck, Dunlop
and Menu (1993). Dunlop, Ferrari and Fontes (2001) proved bounds (slightly
worse than) \reff{eq:b1} for a one dimensional interface related to the phase
separation line in the two dimensional Ising model at zero temperature.
Funaki and Olla (2001) studied a one dimensional model in a finite box
rescaled as the square of the time.

The strategy to show part of Theorem \ref{bounds_gauss} is to compare the wall
process with a ``free process'' --- in our case the serial harness --- as
proposed by Dunlop, Ferrari and Fontes (2001). The following lemmas are the
basic ingredients in this approach.  The first two concern moderate deviations
of the serial harness $Y_n$; they are then extended to the wall process $W_n$
in the last one.

\begin{lm}
\label{dev1}
If the distribution of $\eps$ is in ${\mathcal L}_1^-$, then in $d\leq2$ there
exist constants $k, c, c'>0$ such that for all $K>0$ and $0\leq l\leq n$,
\begin{equation}
    \label{eq:lm2}
\P[Y_l(0)\ge K\sqrt{s(n)\log n}]
\le kn^{c-c'K}.
\end{equation}
\end{lm}

\begin{lm}
\label{dev}
If the distribution of $\eps$ is in
${\mathcal L}_\a^-$ for some $\a\geq1$, then in $d\geq3$
there exist constants $k, c, c'>0$ such that, for all $K>0$
and $0\leq l\leq n$,
\begin{itemize}
\item[(i)] if $\a\ne1+d/2$, then
\begin{equation}
    \label{eq:lm1}
\P[Y_l(0)\ge K(\log n)^{\frac{1}{\a}\vee\frac{2}{2+d}}]\le kn^{c-c'K};
\end{equation}
\item[(ii)] and if $\a=1+d/2$, then
\begin{equation}
    \label{eq:lm1a}
\P[Y_l(0)\ge KL_n(1+2/d)]\le kn^{c-c'K},
\end{equation}
where $L_n(\cdot)$ is defined in~(\ref{eq:qn}) below.
\end{itemize}
\end{lm}

\begin{lm}
\label{devw}
The bounds of Lemmas~\ref{dev} and~\ref{dev1} hold for $l=n$
if we replace $Y_n$ with $W_n$, possibly with worse constants
$k, c$.
\end{lm}

We conclude this introduction with a remark concerning the form~(\ref{I.22})
of the interaction with the wall. Two other choices are also natural.
First, if the noise would push the process below zero, simply do nothing.
Or, in the same case, only take the convex combination without a noise.
Formally, these two cases are, respectively
$$W'_0(i) =W''_0(i)\equiv0,$$
and for $n\geq1$
\begin{eqnarray}
\label{alternative.1}
W'_n(i) = \left\{
            \begin{array}{ll}
        \sum_{j \in \Z^d} p(i,j) W'_{n-1}(j) \, + \, \eps_n(i), & \mbox{if\ this\ is\ positive,}\\{}\hspace{2.3cm}
        W'_{n-1}(i),                       & \mbox{otherwise;} \\
        \end{array}
        \right.
\end{eqnarray}
and
\begin{eqnarray}
\label{alternative.2}
W''_n(i) = \left\{
            \begin{array}{ll}
        \sum_{j \in \Z^d} p(i,j) W''_{n-1}(j) \, + \, \eps_n(i),
                                & \mbox{if\ this\ is\ positive,}\\
        \sum_{j \in \Z^d} p(i,j) W''_{n-1}(j),
                                & \mbox{otherwise.} \\
        \end{array}
        \right.
\end{eqnarray}

Coupling $W, W', W''$ by the same realization of the noise variables, one sees that,
stochastically, both $W' \geq W$ and $W'' \geq W$. This implies immediately
that any lower bound for $\mu_n$ (in particular the ones in this paper)
hold for $\mu'_n:=\E W_n'(0)$ and $\mu''_n:=\E W_n''(0)$ as well.
These dominations also imply immediately the validity of the
results of Theorem~\ref{delocalisation} (a-c) for $W'$ and $W''$.
For the analogue of Theorem~\ref{delocalisation} (d), domination
does not help (it goes in the wrong direction). An argument along
the same lines as the one for $W$ can be made for $W''$
straightforwardly; see paragraphs containing~(\ref{I.22bbb})
and~(\ref{I.22bb}).
Under the assumption that ${\mathcal P}(0,0)>0$, one can also make a
similar argument for $W'$; otherwise, the matter is more delicate,
and we do not have an argument.

As for upper bounds for $\mu'_n, \mu''_n$, the ones we get
for $\mu_n$ also hold for both of them, since the
proof only relies on the free process started at some height $r$ dominating
stochastically the wall process started at the same height, and this holds
for all three choices.


\section{Delocalization}
\label{s2}

In this section we show Theorem~\ref{delocalisation}.  The wall process is
\emph{attractive}, that is,
\begin{equation}
  \label{p10}
  \hbox{if }W\le W' \hbox{ then }({\mathcal P}W + \eps_0)^+ \le ({\mathcal
  P}W' + \eps_0)^+ \qquad a.s.
\end{equation}
coordinatewise, which implies
\begin{equation}
  \label{p101}
  \hbox{if }W_n\le W'_n  \hbox{
  stochastically,  then } W_{n+1}\le W'_{n+1} \hbox{
  stochastically.}
\end{equation}
Since for the process with initial flat surface $0\equiv W_0\le W_1$
\emph{a.s.} this implies \reff{p11}.

Theorem \ref{delocalisation} is a consequence of the following three lemmas.
\begin{lm}
\label{p30}
There is no invariant measure for $(W_n)$ with finite mean.
\end{lm}

\proof Suppose there exists an invariant measure $\nu_o$ with finite mean $m_o$.
Let $I = [-c,c]$ be the support of the distribution of $\eps $. Then there
exists $0<c' < c$ such that $\P[\eps < -c'] > 0$ and, by Markov's inequality,
for any $n$, $\P[\sum_j p_n(0,j) W(j) < 2 m_o] > \frac{1}{2}$, where $p_n$
are the $n$-step transition probabilities.

The preceding implies that the process started from the invariant measure $\nu_o$
reaches the wall at the origin in $n' = 2m_o/{c'}$ steps with strictly positive
probability. This yields a positive drift,
contradicting the assumption.
\qed

\begin{lm}
  \label{p31}
  Every invariant measure for $(W_n)$ dominates stochastically
$$\lim_n \P(W_n\in  \cdot).$$
\end{lm}

\proof Attractiveness \reff{p101} implies that the law of $W_n$ is
stochastically non decreasing and hence converges to a limit. Since the
initial flat configuration is dominated by any other, any invariant measure
dominates stochastically that limit.  \square

Consider the family of processes $((W^{k}_n\,,\, n\ge k),\,k\in \Z)$
defined by
\begin{eqnarray}
\label{I.22c}
W^k_n = \left\{
            \begin{array}{lll}{}\hspace{1cm}0, & \mbox{ if } &n =k, \\
                ({\mathcal P}W^k_{n-1}+\eps_n)^+,& \mbox{ if }  & n \geq k+1, \\
         \end{array}
         \right.
\end{eqnarray}
$ (W^k_n\,,\, n\ge k)$ is the wall process evolving from time $k$ on, having
flat configuration at initial time $k$. It is clear that for $k\ge 0$,
\begin{equation}
  \label{p12}
  W^{-k}_0 \stackrel{d}= W^0_k \,(=W_k) \,.
\end{equation}
Since $0=W^k_k\le W^{k-1}_k$, by attractiveness \reff{p10}, $W^k_n\le
W^{k-1}_n$ for all $n\ge k$, and in particular:
\begin{equation}
  \label{p13}
  W^k_0\le
W^{k-1}_0
\end{equation}
so that $W^{-\infty}_0= \lim_{k\to\infty}W^{-k}_0$ is well defined (but could be
infinity).



\begin{lm}
  \label{p33}
$W^{-\infty}_0$ (and hence $W^{-\infty}_n$ for all $n$) is almost
surely identically infinity.
\end{lm}

\proof
The event $\{W^{-\infty}_0=\infty\}$ belongs to the tail $\s$-algebra of
$\{\eps_k \,:\, k\le 0\}$, and is thus trivial. Write
\begin{eqnarray}
\label{de3}
W^{-\infty}_0&=&(\eps_0+\p W^{-\infty}_{-1})^+=\ldots\\
&=&
(\eps_0+\p(\eps_{-1}+\ldots\p(\eps_{-k+1}+\p W^{-\infty}_{-k})^+\ldots)^+)^+\\
&\geq&U_k+\p^{k}W^{-\infty}_{-k},
\end{eqnarray}
for $k>0$, where $U_k=\sum_{i=0}^{k-1}\p^i\eps_{-i}$.  Notice that $U_k$ is
symmetric and that $U_k$ and $W^{-\infty}_{-k}$ are independent: $U_k$ is a
function of $(\eps_i\,:\,-k+1\le i\le 0)$ while $W^{-\infty}_{-k}$ is function of
$(\eps_i\,:\,i\le -k)$.  Since $W^{-\infty}_{-k} \stackrel{d}=W^{-\infty}_{0}$, for
all $k\geq0$
\begin{equation}
\label{de4}
W^{-\infty}_{0}\geq V_k+\p^{k}W^{-\infty}_{0}, \;\hbox{stochastically}
\end{equation}
with $V_k\stackrel{d}= U_k$, $V_k$ and $W^{-\infty}_{0}$ independent.

A key observation is that $W^{-\infty}_{0}$ is ergodic for spatial shifts. This
follows from the fact that $W^{-\infty}_{0}$ is a function of $\eps_n(i)$'s for a
cone of indices $(n,i)$ in $-\N\times\Z^d$ with vertex in $(0,x)$. Now,
$\E(W^{-\infty}_{0})=\infty$, the Ergodic Theorem implies that
$\p^{k}W^{-\infty}_{0}\to\infty$ almost surely as $k\to\infty$.
Indeed,
\begin{equation}
\label{de4a}
\p^{k}W^{-\infty}_{0}(0)=\sum_{i\in\Z^d}p_k(i)W^{-\infty}_{0}(i)
\geq \frac{c}{k^{d/2}}\sum_{|i|\leq{\sqrt k}}W^{-\infty}_{0}(i)\to\infty,
\end{equation}
as $k\to\infty$, by the Ergodic Theorem. We have used the positivity
of $W^{-\infty}_{0}$ and the well known Local Central Limit Theorem estimate
to the effect
that $\inf_{|i|\leq{\sqrt k}}\,p_k(i)\geq c/k^{d/2}$ for some $c>0$. 
For this estimate,
aperiodicity is required; we leave the necessary and straightforward
adaptations for the periodic case to the reader.

Now,~(\ref{de4}),~(\ref{de4a}) and the symmetry of $V_k$ imply that for
arbitrary $M>0$
\begin{eqnarray}
\label{de5}
\P(W^{-\infty}_{0}>M)&\geq&\liminf_{k\to\infty}\P(V_k+\p^{k}W^{-\infty}_{0}>M)\\
&\geq&
\liminf_{k\to\infty}\P(V_k\geq0)\P(\p^{k}W^{-\infty}_{0}>M)\\
&\geq&
\frac12\liminf_{k\to\infty}\P(\p^{k}W^{-\infty}_{0}>M)=\frac12.
\end{eqnarray}
Thus $\P(W^{-\infty}_{0}=\infty)\geq1/2$ and triviality implies
$\P(W^{-\infty}_{0}=\infty)=1$. \qed

\vskip 2mm \noindent{\bf Proof of Theorem \ref{delocalisation}. }  (a) is
immediate consequence of Lemmas \ref{p31} and \ref{p33}: any invariant surface
dominates stochastically $W^{-\infty}_0$ and $W^{-\infty}_0$ is almost surely
identically infinity. (b) follows from Lemma \ref{p33} and \reff{p12}. (c)
follows from the identity $\mu_n = \E W_n(0) = \E W^{-n}_0$ and
the monotone convergence theorem.  Finally, in \reff{I.22bb} below it is shown
that
 \begin{eqnarray}
   \label{I.22bbb}
   \mu_n -\mu_{n-1} &=&  \E\int_{{\mathcal P}W_{n-1}}^{\infty}   \P(\eps>x)\, dx
 \end{eqnarray}
 Since $\eps$ is integrable and ${\mathcal P}W_{n-1}$ increases to infinity in
 probability, \reff{I.22bbb} converges to zero, and we get (d).  \qed


\section{A generic lower bound}
\label{s3}

{}From \reff{I.22a},
\begin{eqnarray}
  \nonumber
  W_n(i)&=&({\mathcal P}W_{n-1}(i)+\eps_n(i))^+\\
\label{I.22b}
&=&{\mathcal P}W_{n-1}(i)+\eps_n(i)
+(-{\mathcal P}W_{n-1}(i)-\eps_n(i))^+
\end{eqnarray}
 Taking expectations, since $\eps$ is symmetric,
 \begin{eqnarray}
   \label{I.22bb}
   \mu_n &=& \mu_{n-1} + \E\int_{{\mathcal P}W_{n-1}}^{\infty}   \P(\eps>x)\, dx.
 \end{eqnarray}
As 
$\int_y^\infty \P(\eps>x) \,dx$ is a convex function of $y$,
\begin{equation}
\label{eq:mu}
 \mu_n\;\ge\; \mu_{n-1} + \int_{\E({\mathcal P}W_{n-1})}^{\infty}   \P(\eps>x)\,
 dx\; =\;  \mu_{n-1}+\E(\eps -\mu_{n-1})^+.
\end{equation}
For $s\geq0$, let $G(s)=\E(\eps -s)^+$, $H(s)=s+G(s)$, and
$\nu(t)$ be such that
$\int_{0}^{\nu(t)}[G(s)]^{-1}\,ds=t.$

\begin{theo}
\label{mn}
$\mu_n\ge\nu(n)$ for all $n\geq0$.
\end{theo}

\begin{rmk}\rm
  This general lower bound does not depend on the dimension.
\end{rmk}

\begin{corol}
\label{logpol+}
If the distribution of $\eps$ belongs to ${\mathcal L}_\a^+$ for some
$\a>0$,
then there exists $c_2=c_2(\alpha)>0$ such that
\begin{equation}
  \label{eq:lbexp}
 \mu_n\ge  c_2 (\log n)^{\frac{1}{\alpha}}.
\end{equation}
\end{corol}

\begin{corol}
\label{poly}
Suppose that the distribution of $\eps$ decays at most polynomially,
i.e.~$\P(\eps>x)\geq c_0x^{-\a}$ for all $x>1$ and some positive constants
$c_0$ and $\a>1$. Then there exists $c_1=c_1(\alpha)>0$ such that
\begin{equation}
  \label{eq:poly}
 \mu_n\ge c_1n^{\frac{1}{\alpha}}.
\end{equation}
\end{corol}

\vskip 2mm \noindent{\bf Proof of Theorem~\ref{mn}. }
Notice first that $\nu(t)$ is a solution of
$$\nu(t)=\int\limits_{0}^{t}G(\nu(s))\,ds$$ and thus
 satisfies
$$\nu(n)=\nu(n-1)+\int\limits_{n-1}^{n}G(\nu(s))ds.$$

Notice also that $G(x)$ is decreasing and $H(x)$ is
increasing. We prove  the lemma by induction. First,
$\mu_0=\nu(0)=0$. Suppose that $\mu_{n-1}\ge\nu(n-1)$. Then,
\begin{eqnarray*}
\nu(n)&=&\nu(n-1)+\int\limits_{n-1}^{n}G(\nu(s))ds\le \nu(n-1)+G(\nu(n-1))\\
&=&H(\nu(n-1))\le H(\mu_{n-1})\le\mu_n,
\end{eqnarray*}
where the last inequality is~(\ref{eq:mu}). \qed

\vskip 2mm \noindent{\bf Proof of Corollary~\ref{poly}. }
Note that \[G(x)=\E(\eps-x)^+=-\int\limits_{x}^{+\infty}(y-x)d\P[\eps\ge y]=
\int\limits_{x}^{+\infty}\P[\eps\ge y]dy\]
and thus
\begin{equation}
  \label{eq:g}
  g(t):=\int\limits_{0}^{t}\frac{ds}{G(s)}=
\int\limits_{0}^{t}\frac{ds}{\int_{s}^{+\infty}\P[\eps\ge y]dy }.
\end{equation}
Thus, from the assumption in the statement of Corollary~\ref{poly},
\begin{equation}
  \label{eq:g1}
g(t)\leq\frac1{c_0}\int\limits_{0}^{t}\frac{ds}{\frac{1}{\alpha-1}s^{1-\alpha}}=
\frac{\alpha-1}{c_0\alpha}t^{\alpha}
\end{equation}
 and
 \[\nu(t)\geq c_1t^{\frac{1}{\alpha}}\]
follows immediately. \qed

 \medskip

\vskip 2mm \noindent{\bf Proof of Corollary~\ref{logpol+}. }
As above, we have
\begin{equation}
  \label{eq:g2}
g(t)\leq\frac1{c}
\int\limits_{0}^{t}\frac{ds}{\int_{s}^{+\infty}e^{-c'y^{\a}}dy}\leq
c_1
\int\limits_{0}^{t}e^{c_2s^{\a}}ds\leq
c_3e^{c_4t^{\a}}
\end{equation}
and the result follows. \qed

 \medskip


\section{Moderate deviations for the serial harness}
The proofs of Lemmas~\ref{dev1} and~\ref{dev} are based on the behavior of
$\E(e^{\l Y_n(0)})$ for small and large $\l$, established in
Lemmas~\ref{cheb1} and~\ref{gen} below.

\begin{lm}
\label{cheb1}
Let $\ln$ be a sequence
of positive numbers such that
\begin{equation}
  \label{eq:bln}
\bln:=\ln/\sqrt{s(n)}\leq1.
\end{equation}
Then there exists a constant
$c$ such that for all $0\leq l\leq n$
\begin{equation}
    \label{eq:ch1}
\E[e^{\bln Y_l(0)}]\leq e^{c\ln^2}.
\end{equation}
\end{lm}
\proof
For all $0\leq l\leq n$
\begin{eqnarray}
\nonumber
\E\Big[e^{\bln{Y}_l(0)}\Big]&=&\prod_{r=0}^{l-1}\prod_{j\in\Z^d}
\E\Big[e^{\bln p_{r}(j)\eps}\Big]
\leq\prod_{r=0}^{l-1} \prod_{j\in\Z^d} e^{c\bln^2p_{r}(j)^2}\\
\nonumber
&=&\exp\{c\ln^2s(n)^{-1}s(l)\}
\leq e^{c\ln^2},
\end{eqnarray}
where $c=\E(e^{\eps})$ and we have used that for a symmetric random
variable $W$, if $|\lambda|\leq1$, then
\begin{equation}
  \label{eq:sym}
\E(e^{\lambda W})\leq 1+\E(e^{W})\lambda^2\leq e^{\E(e^{W})\lambda^2}
\end{equation}
and the fact that $s(\cdot)$ is nondecreasing.
\qed

\medskip

\vskip 2mm \noindent{\bf Proof of Lemma~\ref{dev1}. }
\begin{eqnarray*}
\P[Y_l(0)\ge K\sqrt{s(n)\log n}]
=\P[\bln Y_n(0)\ge \log n^{c'K}]
\le n^{-c'K}\E[e^{\bln Y_n(0)}],
\end{eqnarray*}
where $\ln=c''\sqrt{\log n}$, for an appropriate constant $c''$, and
Lemma~\ref{cheb1} yields the result.
\qed

\bigskip

For the proof of Lemma~\ref{dev}, we will use that in $d\geq3$
\begin{equation}
  \label{eq:s}
  s:=\lim_{n\to\infty}s(n)<\infty.
\end{equation}
We will also need the following converse of~(\ref{eq:sym}).
\begin{lm}
\label{cheb2}
If the distribution of $W$ is in ${\mathcal L}_\a^-$ for some $\a>1$, then
there exists a constant $c$ such that
\begin{equation}
  \label{eq:sym2}
\E(e^{\lambda W})\leq e^{c\lambda^\beta},
\end{equation}
for all $\l\geq1$, where $\beta=\a/(\a-1)$.
\end{lm}
\proof
We have that
\begin{equation}
  \label{eq:c22}
\E e^{\lambda W}\leq1+c\int_{0}^{\infty}e^{\l x}e^{-c'x^{\alpha}}dx
=1+c_1\int_{0}^{\infty}e^{\tl x}e^{-x^{\alpha}}dx,
\end{equation}
where $\tl=\l/c^{1/\a}$.
Now, we write the integral in~(\ref{eq:c22}) as
\begin{eqnarray*}
\int_{0}^{(2\tl)^{\b-1}}e^{\tl x}dx
+\int_{(2\tl)^{\b-1}}^{\infty}e^{\tl x-x^{\alpha}}dx.
\end{eqnarray*}
The former integral is bounded above by $e^{c'''\lambda^\beta}$.
The latter one is bounded above by a uniform constant. \qed

\bigskip

\begin{lm}
\label{gen}
In $d\geq3$, if the distribution of $\eps$ is in ${\mathcal L}_\a^-$
for some $\a>1$, then there exist a constant $c$ such that for all large $q$
\begin{equation}
  \label{eq:ge1}
\E(e^{qY_n(0)})\leq
\begin{cases}
e^{cq^{\beta\vee(1+2/d)}},&\mbox{ if }\a\ne1+d/2;\\
e^{cq^{1+2/d}\log q},&\mbox{ if }\a=1+d/2,
\end{cases}
\end{equation}
where $\beta=\a/(\a-1)$ as before.
\end{lm}

\proof
\begin{eqnarray}
\nonumber
\E(e^{qY_n(0)})&\leq&\prod_{k=0}^\infty\prod_{x\in\Z^d}\E(e^{qp_k(x)\eps})
\leq\!\!\!\prod_{k,x:qp_k(x)>1}\!\!\!e^{c(qp_k(x))^\b}\!\!\!
    \prod_{k,x:qp_k(x)\leq1}\!\!\!e^{c(qp_k(x))^2}\\
\label{eq:ge2}
&=&\exp\left\{c\left[\sum_{k,x:qp_k(x)>1}\!\!\!(qp_k(x))^\b
+\sum_{k,x:qp_k(x)\leq1}\!\!\!(qp_k(x))^2\right]\right\}.
\end{eqnarray}
We now estimate the expression within square brackets in~(\ref{eq:ge2}).
If $\b\geq2$ or, equivalently, $1<\a\leq2$, then that expression is
bounded above by
\begin{equation}
  \label{eq:ge3}
  q^\b\sum_{k,x}p^2_k(x)=q^\b s.
\end{equation}

For the case $1<\b<2$ (equivalently, $\a>2$), we use the well known
estimate on $p_k:=\sup_{x\in\Z^d}p_k(x)$:
there exists a constant $C$ such that for all $k\geq1$
\begin{equation}
  \label{eq:pk}
  p_k\leq Ck^{-d/2}
\end{equation}
(see e.g. Spitzer (1976))
to conclude that the expression within square brackets in~(\ref{eq:ge2}) is
bounded above by
\begin{equation}
  \label{eq:ge4}
  q^\b\!\!\sum_{k=0}^{{}\quad(Cq)^{2/d}}p^{\b-1}_k
+q^2\!\!\sum_{k=(Cq)^{2/d}}^{\infty}p_k
\leq
  C'q^\b\!\!\sum_{k=1}^{{}\quad(Cq)^{2/d}}k^{-d(\b-1)/2}
+C''q^{1+2/d}
\end{equation}
for some constants $C',C''$. The result follows. \qed

\bigskip

\vskip 2mm \noindent{\bf Proof of Lemma~\ref{dev}. }
Let $Q_n$ be a sequence of positive numbers such that
$Q_n=o(\log n)$ and $q_n=(\log n)/Q_n$. Then
\begin{equation}
  \label{eq:de1}
  \P[Y_l(0)\ge KQ_n]
  \leq\P[q_nY_l(0)\ge K(\log n)]\leq n^{-K}\E(e^{q_nY_l(0)}).
\end{equation}
We can thus use Lemma~\ref{gen} for $q_n$.
Therefore, if $1<\alpha\ne1+d/2$, making
$Q_n=(\log n)^{\frac{1}{\a}\vee\frac{2}{2+d}}$, we have
$q_n=(\log n)^{1-(\frac{1}{\a}\vee\frac{2}{2+d})}=
(\log n)^{\frac{1}{\b}\wedge\frac{d}{d+2}}$ and
thus, from~(\ref{eq:ge1})
\begin{equation}
  \label{eq:de2}
  \P[Y_l(0)\ge K(\log n)^{\frac{1}{\a}\vee\frac{2}{2+d}}]
  \leq n^{c-K}.
\end{equation}

\bigskip

If $\alpha=1+d/2$, we make
$Q_n=L_n(1+2/d)$, and thus
$q_n=(\log n)/L_n(1+2/d)=\ell_n(1+2/d)$. From~(\ref{eq:ge1})
and the definition of $\ell_n(1+2/d)$ (above~(\ref{eq:qn}) below)
\begin{equation}
  \label{eq:de3}
  \P[Y_l(0)\ge KL_n(1+2/d)]\leq n^{c-K}.
\end{equation}

\bigskip

For $\alpha=1$,  we have
\begin{eqnarray}
\E e^{Y_n(0)}=\prod_{k,x}\E e^{p_k(x)\eps}
\leq e^{c\sum_{k,x}p^2_k(x)}=e^{cs},
\end{eqnarray}
where we have used~(\ref{eq:sym}). Thus, we obtain that
\[\P[Y_n(0)>K\log n]\le Cn^{-K}.\]

\qed


\section{Moderate deviations for the wall process}

In this section we show Lemma~\ref{devw}.  Introduce new processes $W^{0,r}_n$
and $Y^{0,r}_n$, which have the same evolution as $W_n$, respectively $Y_n$,
but are started at time zero at height $r \in \N$. That is, $W^{0,r}_0(i) =
Y^{0,r}_0(i) = r$, for all $i \in \Z^d$.

Let
\begin{equation}
\label{an}
  a_n =
\begin{cases}
2 K (\log n)^{\frac{1}{\a}\vee\frac{2}{2+d}},&
\mbox{for the extension of~(\ref{eq:lm1});}\\
2 K L_n(1+2/d),&
\mbox{for the extension of~(\ref{eq:lm1a});}\\
2 K \sqrt{s(n)\log n},&
\mbox{for the extension of~(\ref{eq:lm2})}.
\end{cases}
\end{equation}
 Then,
\begin{eqnarray}
\label{UB.2}
\nonumber
&&\P\left[ W_n(0) \geq a_n \right]
\leq \P\left[W^{0,r}_n(0) \geq a_n \right]\\
\nonumber
&=&\P\left[W^{0,r}_n(0)\geq a_n, W^{0,r}_n(0) = Y^{0,r}_n(0)\right]\\
\nonumber
&+&\P\left[W^{0,r}_n(0) \geq a_n, W^{0,r}_n(0) \neq Y^{0,r}_n(0)\right]\\
\label{eq:ub21}
&\leq&\P\left[Y^{0,r}_n(0) \geq a_n\right]\\
\label{eq:ub22}
&+&\P\left[W^{0,r}_n(0) \neq Y^{0,r}_n(0)\right]
\end{eqnarray}
To get a bound for the probability in~(\ref{eq:ub21}) of the
form~(\ref{eq:lm2}-\ref{eq:lm1a}), we take $r=a_n/2$ and
use~(\ref{eq:lm2}-\ref{eq:lm1a}).

The probability in~(\ref{eq:ub22}) is treated as follows.
Note that $W^{0,r}_n(0)$
and $Y^{0,r}_n(0)$ differ if a discrepancy occurs in the cone ($v$ is the maximal
speed of a discrepancy)
\begin{equation}
\label{UB.7}
\{(l,j) \in \N_0 \times \Z^d: l\leq n, \quad |j| \leq v(n - l)\},
\end{equation}
that is,
\begin{equation}
\nonumber
\{Y^{0,r}_n(0)\neq W^{0,r}_n(0)\}
=\{Y^{0,r}_l(j)<0\mbox{ for some }(l,j)\mbox{ with }
                l\leq n,|j|\leq v(n-l)\}.
\end{equation}
Since $Y^{0,r}_n(0)$ has the same law as $Y_n(0)+r$ and by symmetry, we have
\begin{equation}
\label{UB.9}
\P[Y^{0,r}_l(j) < 0] = \P[Y_l(j) < - r] = \P[Y_l(j) > r].
\end{equation}
Hence,
\begin{eqnarray}
\label{UB.10}
\P[Y^{0,r}_n(0) \leq W^{0,r}_n(0)]
&=&\P[\exists\, (l,j) \mbox{ with } l \leq n, |j| \leq v(n-l): Y_l(j) > r] \nonumber\\
&\leq& \sum_{l=0}^{n} \sum_{|j| \leq v(n-l)}
                \P[Y_l(j) > r]. \nonumber
\end{eqnarray}
Taking $r=a_n/2$ as before and using~(\ref{eq:lm2}-\ref{eq:lm1a}), we obtain
\begin{equation}
\label{UB.10bis}
\P[Y^{0,r}_n(0) \neq W^{0,r}_n(0)]
\leq kn^{c-c'K} \sum_{l=0}^{n} \sum_{|j| \leq v(n-l)}1
\leq  k'n^{c''-c'K},
\end{equation}
for some $k',c''$.
\qed


\section{Bounds for the wall process}

For $\gamma>1$, define $\ell_n(\gamma)$ as the solution of $x^{\gamma}\log
x=\log n$, and let
\begin{equation}
  \label{eq:qn}
L_n(\gamma)=(\log n)/\ell_n(\gamma).
\end{equation}
Note that
\begin{equation}
  \label{p103}
  (\log n)^{1-\frac1\gamma}\leq L_n(\gamma)\leq(\log n)^{1-\frac1\gamma}
(\log\log n)^{\frac1\gamma} \hbox{ for all }n.
\end{equation}

\begin{theo}
\label{upper_marina}
Suppose that the distribution of $\eps$ belongs to ${\mathcal L}_\a^-$
for some $\a\geq1$.
If $d\ge 3$,
then there exists $c_3=c_3(\alpha, d)>0$ such that
\begin{itemize}
\item[(i)] if $1\le\alpha\ne1+\frac{d}{2}$, then
\begin{equation}
    \label{eq:ub1}
    \mu_n\le c_3(\log n)^{\frac{1}{\alpha}\vee\frac{2}{2+d}};
  \end{equation}
\item[(ii)] if $\alpha=1+\frac{d}{2}$, then for all $\delta>0$ we have
\begin{equation}
    \label{eq:ub2}
    \mu_n\le c_3L_n(1+2/d);
  \end{equation}
\end{itemize}
If $d=2$, then there exists $c_3$ such that
\begin{equation}
    \label{eq:ub4}
   \mu_n\le c_3\log n.
  \end{equation}
\end{theo}
\begin{rmk}\rm
\label{rmk:ub2}
{}From~(\ref{eq:ub2}) and \reff{p103}, a slightly weaker alternative
to~(\ref{eq:ub2}) is
\begin{equation}
    \label{eq:ub2p}
    \mu_n\le c_3(\log n)^{\frac{2}{2+d}}(\log\log n)^{\frac{d}{2+d}}.
  \end{equation}
\end{rmk}

We now restrict attention to the class of exponentially decaying noise
distributions. When the noise distribution is in ${\mathcal L}_\a$, $\a\geq1$, the
results in Corollary~\ref{logpol+} and Theorem~\ref{upper_marina} are our best
explicit bounds (to leading order) for $d\geq3$ and $d=2,\,1\leq\a\leq 2$. For
$d=1,\a\geq1$ and $d=2,\,\a>2$, we have better bounds, which we discuss now.

\begin{theo}
\label{bounds_pablo_beat}
If the distribution of $\eps$ is in ${\mathcal  L}_1^-$,
then for $d\leq2$, there exist constants $c, C>0$ such that
\begin{equation}
\label{I.24}
c \sqrt{s(n)}\leq\mu_n \leq C\sqrt{s(n)\log n},
\end{equation}
where $s(n)$ is defined in (\ref{I.23bis}).
In particular
\begin{itemize}
\item[(i)] for $d=1$,
\begin{equation}
\label{I.24a}
c n^{1/4}\leq\mu_n \leq C n^{1/4}\sqrt{\log n};
\end{equation}
\item[(ii)]  and for $d=2$,
\begin{equation}
\label{I.25}
c \sqrt{\log n}\leq\mu_n \leq C\log n
\end{equation}
\end{itemize}
\end{theo}

\begin{rmk}\rm
\label{4mom}
The lower bound in~(\ref{I.24}) actually holds under the weaker assumption
that $\E(\eps^2)<\infty$. See Remark~\ref{4moma} below.
\end{rmk}

We prove first the lower bound~(\ref{I.24}).  The first step is to calculate
the variance of the serial harness, which will give us the proper scaling.
{}From (\ref{I.23}) we get (this is already contained in Hammersley (1965))
$\E\,Y_n(0) = 0$ and $\E\,Y_n(0)^2 =\s^2 s(n)$.

The correct scaling for the serial harness is therefore $s(n)^{1/2}$, and we
define accordingly
\begin{equation}
\label{LB.5}
\tilde{Y}_n(0) \equiv s(n)^{- \frac{1}{2}} Y_n(0).
\end{equation}
Analogously we define $\tilde{W}_n(0)$ for the wall process.  We now show
that $\tilde{Y}_n(0)$ is uniformly integrable (with respect to $n$).

\begin{lm}
\label{lb.2}
The process $(\tilde{Y}_n(0))_n$ satisfies $\sup_n \E( e^{|\tilde{Y}_n(0)|}) <
\infty$.
\end{lm}

\proof
By symmetry of the $\eps$,
$\E(e^{|\tilde{Y}_n(0)|})\leq2\E(e^{\tilde{Y}_n(0)})\leq2e^c$,
where the last inequality follows from Lemma~\ref{cheb1}
with $\ln\equiv1$.
\qed

\medskip

{}From Lemma~(\ref{lb.2}) it follows immediately that $s(n)^{-1}Y_n(0)^2$ is
uniformly integrable.

\begin{lm}
\label{lb.3}
  There exists a constant $c>0$ such that for all $n$
\begin{equation}
\label{LB.9}
\E|\tilde{Y}_n(0)| > c.
\end{equation}
\end{lm}

\proof
Clearly, for any positive $M$,
\begin{eqnarray}
\label{LB.10}
\E[\tilde{Y}_n(0)^2]
&=& \E\left[\tilde{Y}_n(0)^2 \1\{|\tilde{Y}_n(0)| > M\}\right]
        + \E\,\left[\tilde{Y}_n(0)^2 \1\{\{|\tilde{Y}_n(0)| \leq M\}\}\right]
         \nonumber \\
&\leq& \E\,\left[\tilde{Y}_n(0)^2 \1\{|\tilde{Y}_n(0)| > M\}\right]
        + M \E[|\tilde{Y}_n(0)|].
\end{eqnarray}
Since $\tilde{Y}_n(0)^2$ is uniformly integrable, for each
$\delta>0$ we can choose $M>0$ such that
\begin{equation}
\label{LB.11}
\E\,\left[\tilde{Y}_n(0)^2 \1\{|\tilde{Y}_n(0)| > M\}\right]
< \delta,
\end{equation}
uniformly in $n$. Thus
\begin{equation}
\label{LB.12}
\E\,[|\tilde{Y}_n(0)|]
\geq \frac{\E\,[\tilde{Y}_n(0)^2] - \delta}{M} = \frac{\s^2 - \delta}{M} > c > 0,
\end{equation}
for some $\delta > 0$.
\qed

 \medskip

 We finally prove the result about the wall process by coupling it with the
 serial harness using the same disorder variables ${\mathcal E}$.  By symmetry,
\begin{equation}
\label{LB.13}
\E\,[|\tilde{Y}_n(0)|]
= \E\,[(\tilde{Y}_n(0))^+] + \E\,[(-\tilde{Y}_n(0))^+]
= 2 \E\,[(\tilde{Y}_n(0))^+].
\end{equation}
On the other hand, by construction,
$\tilde{W}_n(0) \geq (\tilde{Y}_n(0))^+$, and therefore,
\begin{equation}
\label{LB.14}
\E\,[\tilde{W}_n(0)] \geq \E\,[(\tilde{Y}_n(0))^+]
\geq  \frac{1}{2} \E\,[|\tilde{Y}_n(0)|] \geq  c' > 0.
\end{equation}
This proves the lower bound (\ref{I.24}).

\medskip
The upper bounds~(\ref{eq:ub1}-\ref{eq:ub4}) and~(\ref{I.25})
follow from Lemma~\ref{devw} in the same, following way.
Let $a_n$ be as in~(\ref{an}) and $b_n=a_n/(2K)$.
Then
\begin{eqnarray*}
\mu_n/b_n&=&\E[{W}_n(0)/b_n]=\int_0^\infty\P({W}_n(0)>Kb_n)\,dK\\
&\leq&c/c'+k\int_{c/c'}^\infty n^{c-c'K}\,dK\leq C,
\end{eqnarray*}
for some constant $C$.
\qed

\begin{rmk}\rm
\label{4moma}
The lower bound in~(\ref{I.24}) actually holds under the weaker assumption
that $\E(\eps^2)<\infty$, since this is enough to have
$\tilde{Y}_n(0)^2$ uniformly integrable.
\end{rmk}


\section*{Acknowledgments.}
We thank Servet Mart\'\i nez for many discussions on harnesses, in particular
for pointing out the minimality of $W^{-\infty}_0$. We thank Fran\c cois Dunlop
for references to the physical literature.

This paper is supported by Funda\c c\~ao de Apoio \`a Pesquisa do Estado de
S\~ao Paulo (FAPESP), Conselho Nacional de Desenvolvimento Cient\'{\i}fico e
Tecnol\'ogico (CNPq), Programa N\'ucleos de Excel\^encia (PRONEX).


\parskip0pt

\obeylines
Pablo A.~Ferrari \hfill Luiz R.~G.~Fontes
IME USP \hfill IME USP
Caixa Postal 66281  \hfill Caixa Postal 66281
05311-970 - S\~{a}o Paulo \hfill 05311-970 - S\~{a}o Paulo
BRAZIL \hfill BRAZIL
{\tt pablo@ime.usp.br} \hfill {\tt lrenato@ime.usp.br}
http://www.ime.usp.br/\~{}pablo  \hfill
\vskip 5mm
Beat M.~Niederhauser \hfill Marina Vachkovskaia
IME USP \hfill IMECC UNICAMP
Caixa Postal 66281  \hfill Caixa Postal 6065
05311-970 - S\~{a}o Paulo \hfill 13081-970 Campinas SP
BRAZIL \hfill BRAZIL
{\tt beat@ime.usp.br} \hfill {\tt marinav@ime.unicamp.br}
http://www.ime.usp.br/\~{}beat  \hfill


\end{document}